\documentclass{icm2010}

\usepackage{amscd}

\title[Embedded contact homology]{Embedded contact homology and its
  applications} \author[Michael Hutchings]{Michael Hutchings
  \thanks{Partially supported by NSF grant DMS-0806037} }
\contact[hutching@math.berkeley.edu]{Mathematics Department, 970 Evans
  Hall, University of California, Berkeley CA 94720 USA}

\newtheorem{theorem}{Theorem}[section]
\newtheorem{proposition}[theorem]{Proposition}

\theoremstyle{definition}

\newtheorem*{remark}{Remark}
\newtheorem{question}[theorem]{Question}

\newcommand{\mc}[1]{{\mathcal #1}}

\numberwithin{equation}{section}

\newcommand{\floor}[1]{\left\lfloor #1 \right\rfloor}

\newcommand{\eqdef}{\;{:=}\;}

\renewcommand{\frak}{\mathfrak}

\newcommand{\C}{{\mathbb C}}

\newcommand{\R}{{\mathbb R}}
\newcommand{\N}{{\mathbb N}}
\newcommand{\Z}{{\mathbb Z}}
\newcommand{\op}{\operatorname}

\newcommand{\Spinc}{\op{Spin}^c}

\newcommand{\Ker}{\op{Ker}}

\newcommand{\bpm}{\begin{pmatrix}}
\newcommand{\epm}{\end{pmatrix}}

\renewcommand{\epsilon}{\varepsilon}

\begin{document}

\begin{abstract}
  Embedded contact homology (ECH) is a kind of Floer homology for
  contact three-manifolds.  Taubes has shown that ECH is isomorphic to
  a version of Seiberg-Witten Floer homology (and both are
  conjecturally isomorphic to a version of Heegaard Floer homology).
  This isomorphism allows information to be transferred between
  topology and contact geometry in three dimensions.  In the present
  article we first give an overview of the definition of embedded
  contact homology.  We then outline its applications to
  generalizations of the Weinstein conjecture, the Arnold chord
  conjecture, and obstructions to symplectic embeddings in four
  dimensions.
\end{abstract}



\maketitle

\section{Embedded contact homology}

\subsection{Floer homology of 3-manifolds}
\label{sec:floer}

There are various kinds of Floer theory that one can associate to a
closed oriented $3$-manifold with a spin-c structure.  In this article
we regard a {\em spin-c structure\/} on a closed oriented 3-manifold
$Y$ as an equivalence class of oriented $2$-plane fields on $Y$ (i.e.\
oriented rank $2$ subbundles of the tangent bundle $TY$), where two
oriented $2$-plane fields are considered equivalent if they are
homotopic on the complement of a ball in $Y$.  The set of spin-c
structures on $Y$ is an affine space over $H^2(Y;\Z)$, which we denote
by $\Spinc(Y)$.  A spin-c structure $\frak{s}$ has a well-defined
first Chern class $c_1(\frak{s})\in 2H^2(Y;\Z)$.  A {\em torsion\/}
spin-c structure $\frak{s}$ is one for which $c_1(\frak{s})$ is
torsion in $H^2(Y;\Z)$.

One version of Floer theory for spin-c $3$-manifolds is the {\em
  Seiberg-Witten Floer cohomology\/}, or {\em monopole Floer
  cohomology\/}, as defined by Kronheimer-Mrowka \cite{km}.  There are
two basic variants of this theory, which are different only for
torsion spin-c structures; the variant relevant to our story is
denoted by $\widehat{HM}^*(Y,\frak{s})$.  Very roughly, this is the
homology of a chain complex which is generated by $\R$-invariant
solutions to the Seiberg-Witten equations on $\R\times Y$, and whose
differential counts non-$\R$-invariant solutions to the Seiberg-Witten
equations on $\R\times Y$ which converge to two different
$\R$-invariant solutions as the $\R$ coordinate goes to $+\infty$ or
$-\infty$.  This cohomology is a relatively $\Z/d$-graded $\Z$-module,
where $d$ denotes the divisibility of $c_1(\frak{s})$ in
$H^2(Y;\Z)/\op{Torsion}$.

The Seiberg-Witten Floer cohomology $\widehat{HM}^*(Y,\frak{s})$ is
conjecturally isomorphic to a second kind of Floer theory, the {\em
  Heegaard Floer homology\/} $HF^+_*(-Y,\frak{s})$ defined by
Ozsv\'{a}th-Szab\'{o} \cite{os}.  The latter, roughly speaking, is
defined by taking a Heegaard splitting of $Y$, with Heegaard surface
$\Sigma$ of genus $g$, and setting up a version of Lagrangian Floer
homology in $\op{Sym}^g\Sigma$ for two Lagrangians determined by the
Heegaard splitting.  Although the definitions of Seiberg-Witten Floer
theory and Heegaard Floer theory appear very different, there is
extensive evidence that they are isomorphic, and a program for proving
that they are isomorphic is outlined in \cite{lee}.

Seiberg-Witten and Heegaard Floer homology have had a wealth of
applications to three-dimensional topology.  The present article is
concerned with a third kind of Floer homology, called ``embedded
contact homology'' (ECH), which is defined for contact
3-manifolds.  Because ECH is defined directly in terms of contact
geometry, it is well suited to certain applications in this area.

\subsection{Contact geometry preliminaries}

Let $Y$ be a closed oriented 3-manifold.  A {\em contact form\/} on
$Y$ is a $1$-form $\lambda$ on $Y$ such that $\lambda\wedge
d\lambda>0$ everywhere. The contact form $\lambda$ determines a
2-plane field $\xi=\Ker(\lambda)$, oriented by $d\lambda$; an oriented
$2$-plane field obtained in this way is called a {\em contact
  structure\/}.  The contact form $\lambda$ also determines a vector
field $R$, called the {\em Reeb vector field\/}, characterized by
$d\lambda(R,\cdot)=0$ and $\lambda(R)=1$.

Two basic questions are: First, given a closed oriented 3-manifold
$Y$, what is the classification of contact structures on $Y$ (say, up
to homotopy through contact structures)?  Second, given a contact
structure $\xi$, what can one say about the dynamics of the Reeb
vector field for a contact form $\lambda$ with $\Ker(\lambda)=\xi$?
The first question is a subject of active research which we will not
say much about here, except to note that a fundamental theorem of
Eliashberg \cite{elot} implies that every closed oriented 3-manifold
has a contact structure, in fact a unique ``overtwisted'' contact
structure in every homotopy class of oriented 2-plane fields.  (A
contact structure $\xi$ on a 3-manifold $Y$ is called {\em
  overtwisted\/} if there is an embedded disk $D\subset Y$ such that
$TD|_{\partial D} = \xi|_{\partial D}$.)  For more on this topic see
e.g.\ \cite{geiges}.

To discuss the second question, we need to make some definitions.
Given a closed oriented 3-manifold with a contact form, a {\em Reeb
  orbit\/} is a periodic orbit of the Reeb vector field $R$, i.e.\ a
map $\gamma:\R/T\Z\to Y$ for some $T>0$, such that
$\gamma'(t)=R(\gamma(t))$.  Two Reeb orbits are considered equivalent
if they differ by reparametrization.  If $\gamma:\R/T\Z\to Y$ is a
Reeb orbit and if $k$ is a positive integer, the {\em $k^{th}$
  iterate\/} of $\gamma$ is defined to be the pullback of $\gamma$ to
$\R/kT\Z$.  Every Reeb orbit is either embedded in $Y$, or the
$k^{th}$ iterate of an embedded Reeb orbit for some $k>1$.  Given a
contact structure $\xi$, one can ask: What is the minimum number of
embedded Reeb orbits that a contact form $\lambda$ with
$\Ker(\lambda)=\xi$ can have?  Must there exist Reeb orbits with
particular properties?  Some questions of this nature are discussed in
\S\ref{sec:gw} below.

Continuing with the basic definitions, if $\gamma$ is a Reeb orbit as
above, then the linearization of the Reeb flow near $\gamma$ defines
the {\em linearized return map\/} $P_\gamma$, which is an automorphism
of the two-dimensional symplectic vector space
$(\xi_{\gamma(0)},d\lambda)$.  The Reeb orbit $\gamma$ is called {\em
  nondegenerate\/} if $P_\gamma$ does not have $1$ as an eigenvalue.
If $\gamma$ is nondegenerate, then either $P_\gamma$ has eigenvalues
on the unit circle, in which case $\gamma$ is called {\em elliptic\/};
or else $P_\gamma$ has real eigenvalues, in which case $\gamma$ is
called {\em hyperbolic\/}.  These notions do not depend on the
parametrization of $\gamma$.  We say that the contact form $\lambda$
is {\em nondegenerate\/} if all Reeb orbits are nondegenerate.  For a
given contact structure $\xi$, this property holds for ``generic''
contact forms $\lambda$.

To a nondegenerate Reeb orbit $\gamma$ and a trivialization $\tau$ of
$\gamma^*\xi$, one can associate an integer $\op{CZ}_\tau(\gamma)$
called the {\em Conley-Zehnder index\/}.  Roughly speaking this
measures the rotation of the linearized Reeb flow around $\gamma$ with
respect to $\tau$.  In particular, if $\gamma$ is elliptic, then the
trivialization $\tau$ is homotopic to one with respect to which the
linearized Reeb flow around $\gamma$ rotates by angle $2\pi\theta$ for
some $\theta\in\R\setminus\Z$, and
\[
\op{CZ}_\tau(\gamma) = 2\floor{\theta} + 1,
\]
where $\floor{\cdot}$ denotes the greatest integer function.

\subsection{The ECH chain complex}
\label{sec:ECC}

With the above preliminaries out of the way, we can now define the
embedded contact homology of a closed oriented 3-manifold $Y$ with a
nondegenerate contact form $\lambda$.

To start, define an {\em orbit set\/} to be a finite set of pairs
$\alpha=\{(\alpha_i,m_i)\}$, where the $\alpha_i$'s are distinct
embedded Reeb orbits, and the $m_i$'s are positive integers, which one
can regard as ``multiplicities''.  The orbit set is called {\em
  admissible\/} if $m_i=1$ whenever the Reeb orbit $\alpha_i$ is
hyperbolic.  The homology class of the orbit set $\alpha_i$ is defined
by $[\alpha]\eqdef \sum_im_i\alpha_i\in H_1(Y)$.  Given $\Gamma\in
H_1(Y)$, we define the ECH chain complex $C_*(Y,\lambda,\Gamma)$ to be
the free $\Z$-module generated by admissible orbit sets $\alpha$ with
$[\alpha]=\Gamma$.  As explained in \S\ref{sec:index} below, this chain
complex has a relative $\Z/d$-grading, where $d$ denotes the
divisibility of $c_1(\xi)+2\op{PD}(\Gamma)$ in
$H^2(Y;\Z)/\op{Torsion}$.  We sometimes write a generator $\alpha$ as
above using the multiplicative notation
$\alpha=\prod_i\alpha_i^{m_i}$, although the chain complex grading and
differential that we will define below are not well behaved with
respect to this sort of ``multiplication''.

To define the differential on the chain complex, choose an almost
complex structure $J$ on $\R\times Y$ such that $J$ sends $\partial_s$
to the Reeb vector field $R$, where $s$ denotes the $\R$ coordinate;
$J$ is $\R$-invariant; and $J$ sends the contact structure $\xi$ to
itself, rotating positively with respect to $d\lambda$.  For our
purposes a {\em $J$-holomorphic curve\/} in $\R\times Y$ is a map
$u:\Sigma\to\R\times Y$ where $(\Sigma,j)$ is a punctured compact (not
necessarily connected) Riemann surface, and $du\circ j=J\circ du$.  If
$\gamma$ is a Reeb orbit, a {\em positive end\/} of $u$ at $\gamma$ is
an end of $\Sigma$ on which $u$ is asymptotic to $\R\times\gamma$ as
$s\to+\infty$; a {\em negative end\/} is defined analogously with
$s\to-\infty$.  If $\alpha=\{(\alpha_i,m_i)\}$ and
$\beta=\{(\beta_j,n_j)\}$ are two orbit sets with $[\alpha]=[\beta]\in
H_1(Y)$, let $\mc{M}^J(\alpha,\beta)$ denote the moduli space of
$J$-holomorphic curves as above with positive ends at covers of
$\alpha_i$ with total covering multiplicity $m_i$, negative ends at
covers of $\beta_j$ with total covering multiplicity $n_j$, and no
other ends.  We declare two such $J$-holomorphic curves to be
equivalent if they represent the same {\em current\/} in $\R\times Y$.
For this reason we can identify an element of $\mc{M}^J(\alpha,\beta)$
with the corresponding current in $\R\times Y$, which we typically
denote by $C$.  Note that since $J$ is assumed to be $\R$-invariant,
it follows that $\R$ acts on $\mc{M}^J(\alpha,\beta)$ by translation
of the $\R$ coordinate on $\R\times Y$.

To each $J$-holomorphic curve $C\in\mc{M}^J(\alpha,\beta)$ one can
associate an integer $I(C)$, called the ``ECH index'', which is
explained in \S\ref{sec:index} below.  The differential on the ECH
chain complex counts $J$-holomorphic curves with ECH index $1$, modulo
the $\R$ action by translation.  Curves with ECH index $1$ have
various special properties (assuming that $J$ is generic).  Among
other things, we will see in Proposition~\ref{prop:Rinv} below that if
$I(C)=1$ then $C$ is embedded in $\R\times Y$ (except that $C$ may
contain multiply covered $\R$-invariant cylinders), hence the name
``embedded'' contact homology.  In addition, one can use
Proposition~\ref{prop:Rinv} to show that if $J$ is generic, then the
subset of $\mc{M}^J(\alpha,\beta)$ consisting of $J$-holomorphic
curves $C$ with $I(C)=1$ has finitely many components, each an orbit
of the $\R$ action.

Now fix a generic almost complex
structure $J$.  One then defines the differential $\partial$ on the
ECH chain complex $C_*(Y,\lambda,\Gamma)$ as follows: If $\alpha$ is
an admissible orbit set with $[\alpha]=\Gamma$, then
\[
\partial\alpha \eqdef \sum_\beta \sum_{\left\{C\in\mc{M}^J(\alpha,\beta)/\R\mid
  I(C)=1\right\}}\varepsilon(C)\cdot \beta.
\]
Here the first sum is over admissible orbit sets $\beta$ with
$[\beta]=\Gamma$.  Also $\varepsilon(C)\in\{\pm1\}$ is a sign,
explained in \cite[\S9]{obg2}; the signs depend on some orientation
choices, but the chain complexes for different orientation choices are
canonically isomorphic to each other.  It is shown in \cite{obg1,obg2}
that $\partial^2=0$.  The homology of this chain complex is the {\em
  embedded contact homology\/} $ECH(Y,\lambda,\Gamma)$.

Although the differential $\partial$ depends on the choice of $J$, the
homology of the chain complex does not.  This is a consequence of the
following much stronger theorem of Taubes \cite{e1,e234}, which was
conjectured in \cite{t3}, and which relates ECH to Seiberg-Witten
Floer cohomology.  To state the theorem, observe that the contact
structure $\xi$, being an oriented $2$-plane field, determines a
spin-c structure $\frak{s}_\xi$, see \S\ref{sec:floer}.  We then have:

\begin{theorem}[Taubes]
\label{thm:echswf}
\label{thm:e} There is an isomorphism of relatively $\Z/d$-graded $\Z$-modules
\begin{equation}
\label{eqn:echswf}
ECH_*(Y,\lambda,\Gamma)\simeq
\widehat{HM}^{-*}(Y,\frak{s}_\xi+\op{PD}(\Gamma)).
\end{equation}
\end{theorem}

Here $d$ denotes the divisibility of
\[
c_1(\xi)+2\op{PD}(\Gamma) = c_1(\frak{s}_\xi+\op{PD}(\Gamma))
\]
in $H^2(Y;\Z)/\op{Torsion}$.  Note that both sides of
\eqref{eqn:echswf} are conjecturally isomorphic to the Heegaard Floer
homology $HF^+_*(-Y,\frak{s}_\xi + \op{PD}(\Gamma))$.

\begin{remark}
Both sides of the isomorphism \eqref{eqn:echswf} in fact have absolute
gradings by homotopy classes of oriented $2$-plane fields
\cite{ir,km}, and it is reasonable to conjecture that the isomorphism
\eqref{eqn:echswf} respects these absolute gradings.
\end{remark}

\begin{remark}
  In particular, Theorem~\ref{thm:echswf} implies that, except for
  possible grading shifts, ECH depends only on the 3-manifold $Y$ and
  not on the contact structure.  This is in sharp contrast to the
  symplectic field theory of Eliashberg-Givental-Hofer \cite{egh}
  which, while also defined in terms of Reeb orbits and holomorphic
  curves, is highly sensitive to the contact structure.  In
  particular, the basic versions of symplectic field theory are
  trivial for overtwisted contact structures in three dimensions, see
  \cite{yau,bn}.  On the other hand, while ECH itself does not depend
  on the contact structure, it contains a canonical element which does
  distinguish some contact structures, see \S\ref{sec:empty}.
\end{remark}

\subsection{The ECH index}
\label{sec:index}
To complete the description of the ECH chain complex, we now outline
the definition of the ECH index $I$; full details may be found in
\cite{pfh2,ir}.  This is the subtle part of the definition of ECH, and
we will try to give some idea of its origins.  Meanwhile, on a first
reading one may wish to skip ahead to the examples and applications.

\subsubsection{Four-dimensional motivation}
To motivate the definition of the ECH index, recall that Taubes's
``SW=Gr'' theorem \cite{swgr} relates the Seiberg-Witten invariants of
a closed symplectic 4-manifold $(X,\omega)$, which count solutions to
the Seiberg-Witten equations on $X$, to a ``Gromov invariant'' which
counts certain $J$-holomorphic curves in $X$.  Here $J$ is an
$\omega$-compatible almost complex structure on $X$.  The definition
of ECH is an analogue of Taubes's Gromov invariant for a contact
manifold $(Y,\lambda)$.  Thus Theorem~\ref{thm:echswf} above is an
analogue of SW=Gr for the noncompact symplectic 4-manifold $\R\times
Y$.

For guidance on which $J$-holomorphic curves in $\R\times Y$ to count,
let us recall which $J$-holomorphic curves are counted by Taubes's
Gromov invariant of a closed symplectic 4-manifold $(X,\omega)$.  Let
$C$ be a $J$-holomorphic curve in $(X,\omega)$, and assume that $C$ is
not multiply covered.  If $J$ is generic, then the moduli space of
$J$-holomorphic curves near $C$ is a manifold, whose dimension is a
topological quantity called the {\em Fredholm index\/} of $C$, which
is given by
\begin{equation}
\label{eqn:ind4}
\op{ind}(C) = -\chi(C) + 2 c_1(C).
\end{equation}
Here $c_1(C)$ denotes $\langle c_1(TX),[C]\rangle$, where $TX$ is
regarded as a rank $2$ complex vector bundle via $J$.  In addition, we
have the {\em adjunction formula\/}
\begin{equation}
\label{eqn:adj4}
c_1(C) = \chi(C) + C\cdot C - 2\delta(C).
\end{equation}
Here $C\cdot C$ denotes the self-intersection number of the homology
class $[C]\in H_2(X)$.  In addition $\delta(C)$ is a count of the
singularities of $C$ with positive integer weights, see
\cite[\S7]{mw}, so that $\delta(C)\ge 0$ with equality if and only if
$C$ is embedded.  Now let us define an integer
\begin{equation}
\label{eqn:I4}
I(C) \eqdef c_1(C) + C\cdot C.
\end{equation}
Then equations \eqref{eqn:ind4}, \eqref{eqn:adj4}, and \eqref{eqn:I4}
above imply that
\begin{equation}
\label{eqn:ii4}
\op{ind}(C)\le I(C),
\end{equation}
with equality if and only if $C$ is embedded.  Taubes's Gromov
invariant counts holomorphic currents $C$ with $I(C)=0$, which are
allowed to be multiply covered (but which are not allowed to contain
multiple covers of spheres of negative self-intersection).  Using
\eqref{eqn:ii4}, one can show that if $J$ is generic, then each such
$C$ is a disjoint union of embedded curves of Fredholm index zero,
except that torus components may be multiply covered.  (Multiply
covered tori are counted in a subtle manner explained in \cite{gr}.)

\subsubsection{The three-dimensional story}
\label{sec:index3}

We now consider analogues of the above formulas \eqref{eqn:ind4},
\eqref{eqn:adj4}, and \eqref{eqn:I4} in $\R\times Y$, where
$(Y,\lambda)$ is a contact 3-manifold.  These necessarily include ``boundary
terms'' arising from the ends of the $J$-holomorphic curves.

Let $C\in\mc{M}^J(\alpha,\beta)$ be a $J$-holomorphic curve as in
\S\ref{sec:ECC}, and assume that $C$ is not multiply covered.  It
follows from the main theorem in \cite{dragnev} that if $J$ is
generic, then $\mc{M}^J(\alpha,\beta)$ is a manifold near $C$, whose
dimension can be expressed, similarly to \eqref{eqn:ind4}, as
\begin{equation}
\label{eqn:ind3}
\op{ind}(C) = -\chi(C) + 2c_1(C,\tau) +
\op{CZ}_\tau^0(C).
\end{equation}
Here $c_1(C,\tau)$ denotes the ``relative first Chern class'' of $\xi$
over $C$ with respect to a trivialization $\tau$ of $\xi$ over the
Reeb orbits $\alpha_i$ and $\beta_j$.  This is defined by
algebraically counting the zeroes of a generic section of $\xi$ over
$C$ which on each end is nonvanishing and has winding number zero with
respect to the trivialization $\tau$.  The relative first Chern class
$c_1(C,\tau)$ depends only on $\tau$ and on the relative homology
class of $C$.  Also $\op{CZ}_\tau^0$ denotes the sum, over all the
positive ends of $C$, of the Conley-Zehnder index with respect to
$\tau$ of the corresponding (possibly multiply covered) Reeb orbit,
minus the analogous sum over the negative ends of $C$.

Second, the adjunction formula \eqref{eqn:adj4} is now replaced by the
{\em relative adjunction formula\/}
\begin{equation}
\label{eqn:adj3}
c_1(C,\tau) = \chi(C) + Q_\tau(C) + w_\tau(C) - 2\delta(C).
\end{equation}
Here $Q_\tau(C)$ is a ``relative intersection pairing'' defined in
\cite{pfh2,ir}, which is an analogue of the integer $C\cdot C$ in the
closed case, and which depends only on $\tau$ and the relative
homology class of $C$. Roughly speaking, it is defined by
algebraically counting interior intersections of two generic surfaces
in $[-1,1]\times Y$ with boundary
$\{1\}\times\alpha-\{-1\}\times\beta$ which both represent the
relative homology class of $C$ and which near the boundary have a
special form with respect to the trivialization $\tau$.  As before,
$\delta(C)$ is a count of the singularities of $C$ with positive
integer weights (which is shown in \cite{siefring} to be finite in
this setting).  Finally, $w_\tau(C)$ denotes the {\em asymptotic
  writhe\/} of $C$; to calculate it, take the intersection of $C$ with
$\{s\}\times Y$ where $s>>0$ to obtain a disjoint union of closed
braids around the Reeb orbits $\alpha_i$, use the trivializations
$\tau$ to draw these braids in $\R^3$, and count the crossings with
appropriate signs; then subtract the corresponding count for $s<<0$.

Next we need a new ingredient, which is the following bound on the
asymptotic writhe:
\begin{equation}
\label{eqn:wb}
w_\tau(C) \le \op{CZ}_\tau(\alpha) - \op{CZ}_\tau(\beta) - \op{CZ}_\tau^0(C).
\end{equation}
Here
\[
\op{CZ}_\tau(\alpha) \eqdef
\sum_i\sum_{k=1}^{m_i}\op{CZ}_\tau(\alpha_i^k),
\]
where $\op{CZ}_\tau(\gamma^k)$ denotes the Conley-Zehnder index with
respect to $\tau$ of the $k^{th}$ iterate of $\gamma$.  To prove the
writhe bound \eqref{eqn:wb}, one first needs to understand the
structure of the braids that can arise from the ends of a holomorphic
curve; roughly speaking these are iterated nested cablings of torus
braids, with certain bounds on the winding numbers.  One then needs
some combinatorics to bound the writhes of these braids in terms of
the Conley-Zehnder indices.  The writhe bound was proved in an
analytically simpler situation in \cite{pfh2}; the asymptotic analysis
needed to carry over the proof to the present setting was carried out
by Siefring \cite{siefring}; and an updated proof is given in
\cite{ir}.

Finally, by analogy with \eqref{eqn:I4}, define the {\em ECH index\/}
\begin{equation}
\label{eqn:I3}
I(C) \eqdef c_1(C,\tau) + Q_\tau(C) + \op{CZ}_\tau(\alpha) - \op{CZ}_\tau(\beta).
\end{equation}
One can check that this formula, like the formulas above, does not
depend on the choice of trivialization $\tau$.  It now follows from
\eqref{eqn:ind3}, \eqref{eqn:adj3}, \eqref{eqn:wb}, and \eqref{eqn:I3}
that the {\em index inequality\/}
\begin{equation}
\label{eqn:ii}
\op{ind}(C) \le I(C)
\end{equation}
holds, with equality only if $C$ is embedded.

Recall that we have been assuming in the preceding discussion that $C$ is not
multiply covered.  Without this assumption, one still has the
following proposition, which describes the $I=1$ curves which the ECH
differential counts.

\begin{proposition}
\label{prop:Rinv}
\cite[Cor.\ 11.5]{t3}
Suppose $J$ is generic, and let $C$ be any $J$-holomorphic curve in
$\mc{M}^J(\alpha,\beta)$, possibly multiply covered.  Then:
\begin{description}
\item{(a)}
$I(C)\ge 0$, with equality if and only if $C$ is $\R$-invariant (as a
current).
\item{(b)}
If $I(C)=1$, then $C=C_0\sqcup C_1$, where $I(C_0)=0$, and $C_1$ is
embedded and has $\op{ind}(C_1)=I(C_1)=1$.
\end{description}
\end{proposition}

It may be illuminating to recall the proof here.
As a current, $C$ consists of distinct, irreducible,
non-multiply-covered holomorphic curves $C_1,\ldots,C_k$, covered with
positive integer multiplicities $d_1,\ldots,d_k$.  For simplicity let
us restrict attention to the case when none of the curves $C_i$ is an
$\R$-invariant cylinder.  Let $C'$ be the
holomorphic curve consisting of the union, over $i=1,\ldots,k$, of
$d_i$ different $\R$-translates of $C_i$.  We then have
\begin{equation}
\label{eqn:Rinv}
\sum_{i=1}^kd_i\op{ind}(C_i) = \op{ind}(C') \le I(C')  =
I(C),
\end{equation}
with equality only if the holomorphic curves $C_i$ are embedded and
disjoint.  Here the equality on the left holds because the Fredholom
index is additive under unions, the inequality in the middle is the
index inequality \eqref{eqn:ii} applied to the non-multiply-covered
curve $C'$, and the equality on the right holds because the ECH index
of a holomorphic curve depends only on its relative homology class.
Now since $J$ is generic, and since we made the simplifying assumption
that $C_i$ is not $\R$-invariant, we have $\op{ind}(C_i)>0$ for each
$i$.  We can then read off the conclusions of the proposition in this
case from the inequality \eqref{eqn:Rinv}.

\subsubsection{Grading} The ECH index is also used to define the relative
grading on the ECH chain complex, as follows.  Note that the
definition of the ECH index $I(C)$ depends only on the relative
homology class of $C$, and indeed it makes perfect sense to define
$I(Z)$ as in \eqref{eqn:I3} where $Z$ is any relative homology class
of $2$-chain in $Y$ (not necessarily arising from a $J$-holomorphic
curve) with $\partial Z = \sum_i m_i\alpha_i - \sum_j n_j\beta_j$.  If
$Z'$ is another such relative homology class, then $Z-Z'\in H_2(Y)$,
and one has the {\em index ambiguity formula\/}
\[
I(Z)-I(Z') = \langle c_1(\xi) + 2\op{PD}(\Gamma),Z-Z'\rangle.
\]
We now define the grading difference between two generators $\alpha$
and $\beta$ to be the class of $I(Z)$ in $\Z/d$, where $Z$ is any
relative homology class as above.  The index ambiguity formula shows
that this is well defined, and by definition the differential
decreases the relative grading by $1$.

\subsubsection{Incoming and outgoing partitions and admissibility}
\label{sec:tr}
We now make some technical remarks which will not be needed in the
rest of this article, but which address some frequently asked
questions regarding the definition of ECH.

The first remark is that embeddedness of $C$ is not sufficient for
equality to hold in \eqref{eqn:ii}, unless all of the multiplicities
$m_i$ and $n_j$ equal $1$. A curve $C$ in $\mc{M}^J(\alpha,\beta)$ has
positive ends at covers of $\alpha_i$ with some multiplicities
$q_{i,k}$ whose sum is $\sum_kq_{i,k}=m_i$. If equality holds in
\eqref{eqn:ii}, then the unordered list of multiplicities
$(q_{i,1},q_{i,2},\ldots)$ is uniquely determined by $\alpha_i$ and
$m_i$, and is called the ``outgoing partition''
$P^{\op{out}}_{\alpha_i}(m_i)$.  Likewise the covering multiplicities
associated to the ends of $C$ at covers of $\beta_j$ must comprise a
partition called the ``incoming partition''
$P^{\op{in}}_{\beta_j}(n_j)$.  See e.g.\ \cite[\S4]{ir} for details.
To give the simplest example, if $\gamma$ is an embedded elliptic Reeb
orbit such that the linearized Reeb flow around $\gamma$ with respect
to some trivialization rotates by an angle in the interval $(0,\pi)$,
then $P^{\op{out}}_\gamma(2)=(1,1)$, while
$P^{\op{in}}_{\gamma}(2)=(2)$.

In general, if $\gamma$ is an embedded elliptic Reeb orbit and if $m>1$,
then the incoming and outgoing partitions $P^{\op{in}}_\gamma(m)$ and
$P^{\op{out}}_\gamma(m)$ are always different.  This fact makes the
proof that $\partial^2=0$ quite nontrivial.

On the other hand, suppose $\gamma$ is a hyperbolic embedded Reeb
orbit.  If the linearized return map has positive eigenvalues then
\begin{equation}
\label{eqn:ph}
P^{\op{in}}_\gamma(m) = P^{\op{out}}_\gamma(m)=(1,\ldots,1).
\end{equation}
If the linearized return map has negative eigenvalues then
\begin{equation}
\label{eqn:nh}
P^{\op{in}}_\gamma(m)=P^{\op{out}}_\gamma(m)=\left\{\begin{array}{cl}
    (2,\ldots,2), & \mbox{$m$ even},\\
(2,\ldots,2,1), & \mbox{$m$ odd}.
\end{array}
\right.
\end{equation}
This is one reason why the generators of the ECH chain complex in
\S\ref{sec:ECC} are required to be {\em admissible\/} orbit sets: one
can show using \eqref{eqn:ph} and \eqref{eqn:nh} that if one tries to
glue two $I=1$ holomorphic curves along an inadmissible orbit set,
then there are an even number of ways to glue, which by \cite{bm}
count with cancelling signs.  Thus one must disallow inadmissible
orbit sets in order to obtain $\partial^2=0$.  A similar issue arises
in the definition of symplectic field theory \cite{egh}, where ``bad''
Reeb orbits must be discarded.

\subsection{Example: the ECH of an ellipsoid}
\label{sec:ellipsoid}

We now illustrate the above definitions with what is probably the
simplest example of ECH.  Consider $\C^2=\R^4$ with coordinates
$z_j=x_j+iy_j$ for $j=1,2$.  Let $a,b$ be positive real numbers with
$a/b$ irrational, and consider the ellipsoid
\begin{equation}
\label{eqn:ellipsoid}
E(a,b) \eqdef \left\{(z_1,z_2)\in\C^2\;\bigg|\;\frac{\pi|z_1|^2}{a} +
  \frac{\pi|z_2|^2}{b} \le 1\right\}.
\end{equation}
We now compute the embedded contact homology of $Y=\partial E(a,b)$,
with the contact form
\begin{equation}
\label{eqn:lambda}
\lambda\eqdef\frac{1}{2}\sum_{j=1}^2 (x_jdy_j-y_jdx_j)
\end{equation}
(and of course with $\Gamma=0$).

The Reeb vector field on $Y$ is given by
\[
R = \frac{2\pi}{a}\frac{\partial}{\partial\theta_1} +
\frac{2\pi}{b}\frac{\partial}{\partial\theta_2}
\]
where $\partial/\partial\theta_j\eqdef x_j\partial_{y_j} -
y_j\partial_{x_j}$.  Since $a/b$ is irrational, it follows that there
are just two embedded Reeb orbits $\gamma_1$ and $\gamma_2$, given by
the circles where $z_2=0$ and $z_1=0$ respectively.  These Reeb
orbits, as well as their iterates, are nondegenerate and elliptic.
Indeed there is a natural trivialization $\tau$ of $\xi$ over each
$\gamma_i$ induced by an embedded disk bounded by $\gamma_i$.  With
respect to this trivialization, the linearized Reeb flow around
$\gamma_1$ is rotation by angle $2\pi a/b$, while the linearized Reeb
flow around $\gamma_2$ is rotation by angle $2\pi b/a$.

The generators of the ECH chain complex have the form
$\alpha=\gamma_1^{m_1}\gamma_2^{m_2}$ where $m_1,m_2$ are nonnegative
integers.  We now compute the grading.  The relative $\Z$-grading has
a distinguished refinement to an absolute grading in which the empty
set of Reeb orbits (given by $m_1=m_2=0$ above) has grading $0$.  An
arbitrary generator $\alpha$ as above then has grading
\[
I(\alpha) = c_1(\alpha,\tau) + Q_\tau(\alpha) + \op{CZ}_\tau(\alpha),
\]
where $c_1(\alpha,\tau)$ denotes the relative first Chern class of
$\xi$ over a surface bounded by $\alpha$, and $Q_\tau(\alpha)$ denotes
the relative intersection pairing of such a surface.  Computing using
the above trivialization $\tau$, one finds, see \cite[\S4.2]{wh}, that
\[
\begin{split}
c_1(\alpha,\tau) &= m_1+m_2,\\
Q_\tau(\alpha) &= 2m_1m_2,\\
\op{CZ}_\tau(\alpha) &= \sum_{k=1}^{m_1}(2\floor{ka/b}+1) +
\sum_{k=1}^{m_2}(2\floor{kb/a}+1).
\end{split}
\]
Therefore
\begin{equation}
\label{eqn:lattice}
I(\alpha) = 2\left(m_1+m_2+m_1m_2 + \sum_{k=1}^{m_1}\floor{ka/b} +
  \sum_{k=1}^{m_2}\floor{kb/a}\right).
\end{equation}
In particular, all generators have even grading, so the
differential vanishes, and to determine the homology we just have to
count the number of generators with each grading.

Now if the ECH of $\partial E(a,b)$ is to agree with
$\widehat{HM}^{-*}$ and $HF_*^+$ of $S^3$, then we should get
\begin{equation}
\label{eqn:agree}
ECH_*(\partial E(a,b),\lambda,0) \simeq \left\{\begin{array}{cl} \Z, &
    *=0,2,4,\ldots,
\\ 0, & \mbox{otherwise}.
\end{array}\right.
\end{equation}
It is perhaps not immediately obvious how to deduce this from
\eqref{eqn:lattice}.  The trick is to interpret the right hand side of
\eqref{eqn:lattice} as a count of lattice points, as follows.  Let $T$
denote the triangle in $\R^2$ bounded by the coordinate axes, together
with the line $L$ through the point $(m_1,m_2)$ of slope $-a/b$.  Then
we observe that
\[
I(\alpha) = 2\left(\left|T\cap\Z^2\right|-1\right).
\]
Now if one moves the line $L$ up and to the right, keeping its slope
fixed, then one hits all of the lattice points in the nonnegative
quadrant in succession, each time increasing the number of lattice
points in the triangle $T$ by $1$.  It follows that the ECH chain
complex has one generator in each nonnegative even grading, so
\eqref{eqn:agree} holds.

Usually direct calculations of ECH are not so easy because there are
more Reeb orbits, and one has to understand the holomorphic curves.
But for certain simple contact manifolds this is possible; for example
the ECH of standard contact forms on $T^3$ is computed in \cite{t3},
and these calculations are generalized to $T^2$-bundles over $S^1$ in
\cite{lebow}.

\begin{remark}
  For some mysterious reason, lattice point counts such as the one in
  equation \eqref{eqn:lattice} arise repeatedly in ECH in different
  contexts.  For example one lattice point count comes up in the
  combinatorial part of the proof of the writhe bound \eqref{eqn:wb},
  and in determining the ``partition conditions'' in
  \S\ref{sec:tr}, see \cite[\S4.6]{ir}.  Another lattice
  point count appears in the combinatorial description of the ECH
  chain complex for $T^3$ in \cite[\S1.3]{t3}.
\end{remark}

\subsection{Some additional structures on ECH} 

ECH has various additional structures on it.  We now describe those
structures that are relevant elsewhere in this article.

\subsubsection{The $U$ map}
On the ECH chain complex there is a degree $-2$ chain map
\[
U: C_*(Y,\lambda,\Gamma) \longrightarrow C_{*-2}(Y,\lambda,\Gamma),
\]
see e.g.\ \cite[\S2.5]{wh}.  This is defined similarly to the
differential $\partial$, but instead of counting $I=1$ curves modulo
translation, one counts $I=2$ curves that are required to pass through
a fixed, generic point $z\in\R\times Y$.  This induces a well-defined
map on homology
\[
U: ECH_*(Y,\lambda,\Gamma) \longrightarrow
ECH_{*-2}(Y,\lambda,\Gamma).
\]
Taubes \cite{e5} has shown that this map agrees with an analogous map
on $\widehat{HM}^{-*}$, and it conjecturally agrees with the $U$ map
on $HF^+_*$.  The $U$ map plays a crucial role in the applications to
generalizations of the Weinstein conjecture discussed in
\S\ref{sec:gw} below.

\subsubsection{Filtered ECH}
\label{sec:filtered}
 If $\alpha=\{(\alpha_i,m_i)\}$ is a
generator of the ECH chain complex, define its {\em symplectic
  action\/}
\[
\mc{A}(\alpha) \eqdef \sum_im_i\int_{\alpha_i}\lambda.
\]
It follows from Stokes's theorem and the conditions on the almost
complex structure $J$ that the differential $\partial$ decreases the
symplectic action, i.e.\ if $\langle \partial\alpha,\beta\rangle\neq 0$
then $\mc{A}(\alpha)>\mc{A}(\beta)$.  Given $L\in\R$, we then define
$ECH^L(Y,\lambda,\Gamma)$ to be the homology of the subcomplex of
$C_*(Y,\lambda,\Gamma)$ spanned by generators with symplectic action
less than $L$.  We call this {\em filtered ECH\/}; it is shown in
\cite{cc} that this does not depend on the choice of almost complex
structure $J$.  However, unlike the usual ECH, filtered ECH is not
invariant under deformation of the contact form; see \S\ref{sec:eo}
for some examples.  Filtered ECH has no obvious direct counterpart in
Seiberg-Witten or Heegaard Floer homology, but it plays an important
role in the applications in \S\ref{sec:cc} and \S\ref{sec:eo} below.

\subsubsection{Cobordism maps}
\label{sec:cobordism}
Let $(Y_+,\lambda_+)$ and $(Y_-,\lambda_-)$ be closed oriented
3-manifolds with nondegenerate contact forms. An {\em exact symplectic cobordism\/} from
$(Y_+,\lambda_+)$ to $(Y_-,\lambda_-)$ is a compact symplectic
4-manifold $(X,\omega)$ with boundary $\partial X = Y_+ - Y_-$, such
that there exists a $1$-form $\lambda$ on $X$ with $d\lambda=\omega$
on $X$ and $\lambda|_{Y_\pm}=\lambda_\pm$.  It is shown in \cite{cc}
that an exact symplectic cobordism as above induces maps on filtered
ECH,
\[
\Phi^L(X,\omega): ECH^L(Y_+,\lambda_+;\Z/2) \longrightarrow
ECH^L(Y_-,\lambda_-;\Z/2),
\]
satisfying various axioms.  Here $ECH(Y_\pm,\lambda_\pm;\Z/2)$ denotes
the ECH with $\Z/2$ coefficients, summed over all $\Gamma\in H_1(Y)$,
and regarded as an ungraded $\Z/2$-module.  One axiom is that if
$L<L'$ then the diagram
\[
\begin{CD}
ECH^L(Y_+,\lambda_+;\Z/2) @>{\Phi^L(X,\omega)}>>
ECH^L(Y_-,\lambda_-;\Z/2)\\
@VVV @VVV\\
ECH^{L'}(Y_+,\lambda_+;\Z/2) @>{\Phi^{L'}(X,\omega)}>>
ECH^{L'}(Y_-,\lambda_-;\Z/2)
\end{CD}
\]
commutes, where the vertical arrows are induced by inclusion of chain
complexes.  Thus the direct
limit
\begin{equation}
\label{eqn:dl}
\Phi(X,\omega)\eqdef \lim_{L\to\infty}\Phi^L(X,\omega):
ECH(Y_+,\lambda_+;\Z/2) \longrightarrow ECH(Y_-,\lambda_-;\Z/2)
\end{equation}
is well-defined.  Another axiom is that this direct limit agrees with
the map $\widehat{HM}^{*}(Y_+;\Z/2)\to \widehat{HM}^*(Y_-;\Z/2)$ on
Seiberg-Witten Floer cohomology induced by $X$, under the isomorphism
\eqref{eqn:echswf}.  Here we are considering Seiberg-Witten Floer
cohomology with $\Z/2$ coefficients, summed over all spin-c structures.

\begin{remark}
  The cobordism maps $\Phi^L(X,\omega)$ are defined in \cite{cc} using
  Seiberg-Witten theory and parts of the isomorphism
  \eqref{eqn:echswf}.  It would be natural to try to give an
  alternate, more direct definition of the cobordism maps
  $\Phi^L(X,\omega)$ by counting $I=0$ holomorphic curves in the
  ``completion'' of $X$ obtained by attaching symplectization ends.
  Note that by Stokes's theorem and the exactness of the cobordism,
  such a map would automatically respect the symplectic action
  filtrations.  However there are technical difficulties with defining
  cobordism maps this way, because the compactified $I=0$ moduli
  spaces include broken holomorphic curves which contain multiply
  covered components with negative ECH index.  (There is no analogue
  in this setting of Proposition~\ref{prop:Rinv}, whose proof made
  essential use of the $\R$-invariance of $J$.)  Examples show that
  such broken curves must sometimes make contributions to the
  cobordism map, but it is not known how to define the contribution in
  general.  Fortunately, the Seiberg-Witten definition of
  $\Phi^L(X,\omega)$ is sufficient for the applications considered
  here.
\end{remark}
 
\subsubsection{The contact invariant}
\label{sec:empty}
The {\em empty set\/} is a legitimate generator of the ECH chain
complex.  By the discussion in \S\ref{sec:filtered} it is a cycle, and
we denote its homology class by
\[
c(\xi) \in ECH_0(Y,\lambda,0).
\]
This depends only on the contact structure, although not just on the
3-manifold $Y$.  Indeed the cobordism maps in \S\ref{sec:cobordism}
can be used to show that $c(\xi)$ is nonzero if there is an exact
symplectic cobordism from $(Y,\xi)$ to the empty set.  On the other
hand the argument in the appendix to \cite{yau} implies that
$c(\xi)=0$ if $\xi$ is overtwisted.  Some new families of contact
3-manifolds with vanishing ECH contact invariant are introduced by
Wendl \cite{wendl}.  It is shown by Taubes \cite{e5} that $c(\xi)$
agrees with an analogous contact invariant in Seiberg-Witten Floer
cohomology, and both conjecturally agree with the contact invariant in
Heegaard Floer homology \cite{oscontact}.

\subsection{Analogues of ECH in other contexts}

One can also define a version of ECH for sutured
3-manifolds with contact structures adapted to the sutures, see
\cite{cghh}.  This conjecturally agrees with the sutured Floer
homology of Juh\'{a}sz \cite{juhasz} and with the sutured version of
Seiberg-Witten Floer homology defined by Kronheimer-Mrowka \cite{km}.

There is also an analogue of ECH, called ``periodic Floer homology'',
for mapping tori of area-preserving surface diffeomorphisms, see e.g.\
\cite{pfh3,lt}.

We remark that no analogue of ECH is currently known for contact
manifolds of dimension greater than three.  In higher dimensions one
expects that if $J$ is generic then all non-multiply-covered
$J$-holomorphic curves are embedded, see \cite{oz}.  In addition no
good analogue of Seiberg-Witten theory is known in higher dimensions.

\section{Applications}

Currently all applications of ECH make use of Taubes's isomorphism
\eqref{eqn:echswf}, together with known properties of Seiberg-Witten
Floer homology, to deduce certain properties of ECH which then have
implications for contact geometry.  It is an interesting open problem
to establish the relevant properties of ECH without using
Seiberg-Witten theory.

\subsection{Generalizations of the Weinstein conjecture}
\label{sec:gw}

The {\em Weinstein conjecture\/} in three dimensions asserts that for
any contact form $\lambda$ on a closed oriented 3-manifold $Y$, there
exists a Reeb orbit.  Many cases of this were proved by various
authors, see e.g.\ \cite{h,ach,ch}, and the general case was proved by
Taubes \cite{tw}.  Indeed the three-dimensional Weinstein conjecture
follows immediately from the isomorphism \eqref{eqn:echswf}, together
with a theorem of Kronheimer-Mrowka \cite{km} asserting that
$\widehat{HM}^*$ is always infinitely generated for torsion spin-c
structures.  The reason is that if there were no Reeb orbit, then the
ECH would have just one generator: the empty set of Reeb orbits.
However to prove the Weinstein conjecture one does not need to use the
full force of the isomorphism \eqref{eqn:echswf}; one just needs a way
of passing from Seiberg-Witten Floer generators to ECH generators,
which is what \cite{tw} establishes.

In \cite{wh} we make heavier use of the isomorphism \eqref{eqn:echswf}
to prove some stronger results.  For example:

\begin{theorem}
\label{thm:elliptic}
  Let $\lambda$ be a nondegenerate contact form on a closed oriented
  connected 3-manifold $Y$ such that all Reeb orbits are elliptic.
  Then there are exactly two embedded Reeb orbits, $Y$ is a sphere or
  a lens space, and the two embedded Reeb orbits are the core circles
  of a genus $1$ Heegaard splitting of $Y$.
\end{theorem}

The idea of the proof is as follows.  Since all Reeb orbits are
elliptic, a general property of the ECH index \cite[Prop.\
1.6(c)]{pfh2} implies that all ECH generators have even grading, so
the ECH differential vanishes.  Since $\widehat{HM}^*$ is nonvanishing
for only finitely many spin-c structures, it follows that all Reeb
orbits represent torsion homology classes.  Estimating the number of
ECH generators in a given index range then shows that there are
exactly two embedded Reeb orbits; otherwise there would be either too
few or too many generators to be consistent with the growth rate of
$\widehat{HM}^*$.  Next, known properties of $\widehat{HM}^*$ imply
that the $U$ map is an isomorphism when the grading is sufficiently
large.  This provides a large supply of $I=2$ holomorphic curves in
$\R\times Y$.  By careful use of the adjunction formula
\eqref{eqn:adj3} one can show that at least one of these holomorphic
curves includes a non-$\R$-invariant holomorphic cylinder.  By
adapting ideas from \cite{hwz3}, one can show that this holomorphic
cylinder projects to an embedded surface in $Y$ which generates a
foliation by cylinders of the complement in $Y$ of the Reeb orbits.
This foliation then gives rise to the desired Heegaard splitting.

Theorem~\ref{thm:elliptic} is used in \cite{wh} to extend the
Weinstein conjecture to ``stable Hamiltonian structures'' (a certain
generalization of contact forms) on 3-manifolds that are not torus
bundles over $S^1$.

In addition, a slight refinement of the proof of
Theorem~\ref{thm:elliptic} in \cite{wh} establishes:

\begin{theorem}
Let $Y$ be a closed oriented 3-manifold with a nondegenerate contact
form $\lambda$.  If $Y$ is not a sphere or a lens space, then there
are at least $3$ embedded Reeb orbits.
\end{theorem}

In fact, examples of contact forms with only finitely many embedded
Reeb orbits are hard to come by, and to our knowledge the following
question is open:

\begin{question}
Is there any example of a contact form on a closed connected oriented 3-manifold
with only finitely many embedded Reeb orbits, other than contact forms
on $S^3$ and lens spaces with exactly two embedded Reeb orbits?
\end{question}

It is shown in \cite{hwz} that for a large class of contact forms on
$S^3$ there are either two or infinitely many embedded Reeb orbits.
It is shown in \cite{ch}, using linearized contact homology, that many
contact structures on 3-manifolds (namely those supported by an open
book decomposition with pseudo-Anosov monodromy satisfying a certain
inequality) have the property that for any contact form, there are
infinitely many free homotopy classes of loops that must contain an
embedded Reeb orbit.

\subsection{The Arnold chord conjecture}
\label{sec:cc}

A {\em Legendrian knot\/} in a contact 3-manifold $(Y,\lambda)$ is a
knot $K\subset Y$ such that $TK\subset \xi|_K$.  A {\em Reeb chord\/}
of $K$ is a Reeb trajectory starting and ending on $K$, i.e.\ a path
$\gamma:[0,T]\to Y$ for some $T>0$ such that $\gamma'(t)=R(\gamma(t))$
and $\gamma(0),\gamma(T)\in K$.  The following theorem, proved in
\cite{cc}, is a version of the Arnold chord conjecture \cite{arnold}.
(This was previously known in some cases from \cite{abbas,mohnke}.)

\begin{theorem}
Let $Y_0$ be a closed oriented 3-manifold with a contact form $\lambda_0$,
and let $K$ be a Legendrian knot in $(Y_0,\lambda_0)$.  Then $K$ has a
Reeb chord.
\end{theorem}

The idea of the proof is follows.  Following Weinstein
\cite{weinstein}, one can perform a ``Legendrian surgery'' along $K$
to obtain a new contact manifold $(Y_1,\lambda_1)$, together with an
exact symplectic cobordism $(X,\omega)$ from $(Y_1,\lambda_1)$ to
$(Y_0,\lambda_0)$.  If $K$ has no Reeb chord, and if $\lambda_0$ is
nondegenerate, then one can carry out the Legendrian surgery
construction so that $\lambda_1$ is nondegenerate and, up to a given
action, the Reeb orbits of $\lambda_1$ are the same as those of
$\lambda_0$.  Using this observation, one can show that if $K$ has no
Reeb chord and if $\lambda_0$ is nondegenerate, then the cobordism map
\begin{equation}
\label{eqn:ni}
\Phi(X,\omega): ECH(Y_1,\lambda_1;\Z/2) \longrightarrow
ECH(Y_0,\lambda_0;\Z/2)
\end{equation}
from \S\ref{sec:cobordism} is an isomorphism.  Note that this is what
one would expect by analogy with a very special case of the work of
Bourgeois-Ekholm-Eliashberg \cite{bee}, which studies the behavior of
linearized contact homology under Legendrian surgery, possibly in the
presence of Reeb chords.

But the map
\eqref{eqn:ni} cannot be an isomorphism.  The reason is that as
shown in \cite{km}, the corresponding map on Seiberg-Witten Floer
cohomology fits into an exact triangle
\[
\cdots\to
\widehat{HM}^*(Y_0;\Z/2) \to \widehat{HM}^*(Y_1;\Z/2)\to
\widehat{HM}^*(Y_2;\Z/2) \to \widehat{HM}^*(Y_0;\Z/2) \to \cdots
\]
where $Y_2$ is obtained from $Y_0$ by a different surgery along $K$.
However, as noted before, Kronheimer-Mrowka showed that
$\widehat{HM}^*(Y_2;\Z/2)$ is infinitely generated.  This
contradiction proves the chord conjecture when $\lambda_0$ is nondegenerate.

To deal with the case where $\lambda_0$ is degenerate, one can use
filtered ECH to show that in the nondegenerate case, there exists a
Reeb chord with an upper bound on the length, in terms of a
quantitative measure of the failure of the map \eqref{eqn:ni} to be an
isomorphism.  For example, if $\lambda_0$ is nondegenerate and if
\eqref{eqn:ni} is not surjective, then there exists a Reeb chord of
action at most $A$, where $A$ is the infimum over $L\in\R$ such that
the image of $ECH^L(Y_0,\lambda_0;\Z/2)$ in $ECH(Y_0,\lambda_0;\Z/2)$
is not contained in the image of the map \eqref{eqn:ni}.  One can show
that this upper bound on the length of a Reeb chord is suitably
``continuous'' as one changes the contact form.  A compactness
argument then finds a Reeb chord in the degenerate case.

\subsection{Obstructions to symplectic embeddings}
\label{sec:eo}

ECH also gives obstructions to symplectically embedding one compact
symplectic 4-manifold with boundary into another.  We now explain how
this works in the case of ellipsoids as in \eqref{eqn:ellipsoid}, with
the standard symplectic form $\omega=\sum_{j=1}^2dx_jdy_j$ on $\R^4$.

Given positive real numbers $a,b$, and given a positive integer $k$,
define $(a,b)_k$ to be the $k^{th}$ smallest entry in the array
$(ma+nb)_{m,n\in\N}$.  Here in the definition of ``$k^{th}$ smallest''
we count with repetitions.  For example if $a=b$ then
\[
((a,a)_1,(a,a)_2,\ldots) = (0,a,a,2a,2a,2a,3a,3a,3a,3a,\ldots).
\]
We then have:

\begin{theorem}
\label{thm:eo}
If there is a symplectic embedding of $E(a,b)$ into $E(c,d)$, then
\begin{equation}
\label{eqn:eo}
(a,b)_k \le (c,d)_k
\end{equation}
for all positive integers $k$.
\end{theorem}

To prove this, one can assume without loss of generality that $a/b$
and $c/d$ are irrational and that there is a symplectic embedding
$\varphi:E(a,b) \to \op{int}(E(c,d))$.  Now consider the 4-manifold
$X=E(c,d)\setminus\op{int}(\varphi(E(a,b)))$.  One can show that $X$
defines an exact symplectic cobordism from $\partial E(c,d)$ to
$\partial E(a,b)$, where the latter two 3-manifolds are endowed
with the contact form \eqref{eqn:lambda}.  Since $X$ is diffeomorphic
to the product $[0,1]\times S^3$, the induced map from the
Seiberg-Witten Floer cohomology of $\partial E(c,d)$ to that of
$\partial E(a,b)$ must be an isomorphism.  Recall from \eqref{eqn:dl}
that this map is the direct limit of maps on filtered ECH.  Since the
ECH differentials vanish, it follows that for each $L\in\R$, the
number of ECH generators of $\partial E(c,d)$ with action less than
$L$ does not exceed the number of ECH generators of $\partial E(a,b)$
with action less than $L$.  Since the embedded Reeb orbits in
$\partial E(a,b)$ have action $a$ and $b$, and the embedded Reeb
orbits in $\partial E(c,d)$ have action $c$ and $d$, it follows that
\begin{equation}
\label{eqn:eo2}
\left|\left\{(m,n)\in\N^2\mid cm+dn<L\right\}\right| \le
\left|\left\{(m,n)\in\N^2\mid am+bn < L\right\}\right|.
\end{equation}
The statement that the above inequality holds for all $L\in\R$ is
equivalent to \eqref{eqn:eo}.

For example, if $L$ is large with respect to $a,b,c,d$, then the
inequality \eqref{eqn:eo2} implies that
\[
\frac{L^2}{2cd} \le \frac{L^2}{2ab} + O(L).
\]
We conclude that $ab\le cd$, which is simply the condition that the
volume of $E(a,b)$ is less than or equal to the volume of $E(c,d)$,
which of course is necessary for the existence of a symplectic embedding.  But
taking suitable small $L$ often gives stronger conditions.

The amazing fact is that, at least for the problem of embedding
ellipsoids into balls, the obstruction in Theorem~\ref{thm:eo} is
sharp.  Namely, for each positive real number $a$, define $f(a)$ to be
the infimum over all $c\in\R$ such that $E(a,1)$ symplectically embeds
into the ball $E(c,c)$.  It follows from Theorem~\ref{thm:eo} that
\begin{equation}
\label{eqn:f}
f(a) \ge \sup_{k=2,3,\ldots}\frac{(a,1)_k}{(1,1)_k}.
\end{equation}
On the other hand, McDuff-Schlenk \cite{ms} computed the function $f$
explicitly, obtaining a complicated answer involving Fibonacci
numbers.  Using the result of this calculation, they checked that the
opposite inequality in \eqref{eqn:f} holds.

\begin{question}
  Is there a direct explanation for this?  Does this generalize?  For
  example, does $E(a,b)$ symplectically embed into
  $E(c+\epsilon,d+\epsilon)$ for all $\epsilon>0$ if
  $(a,b)_k\le(c,d)_k$ for all positive integers $k$?
\end{question}

By more involved calculations, one can use ECH to find explicit (but
subtle, number-theoretic) obstructions to symplectic embeddings
involving other simple shapes such as four-dimensional polydisks.  A
systematic treatment of the symplectic embedding obstructions arising
from ECH is given in \cite{qech}.

\end{document}